\documentclass[12pt,bezier]{article}
\usepackage{amsmath,amssymb,amsfonts}

\newcommand{\D}{\displaystyle}
\newcommand{\rg}{\rightarrow}
\newcommand{\st}{\stackrel}

\newcommand{\var}{\varepsilon}

\newcommand{\R}{\Bbb{R}}
\begin{document}
\begin{center}
{\bf \large Analysis of a Disease Transmission Model \\ with two
Groups of Infectives}

\vspace*{1cm}

{\bf M.R. Razvan}\\

{\small\it
Institute for Studies in Theoretical Physics and Mathematics \\
P.O.Box 11395-5746, Tehran, IRAN} \\
E-mail: razvan@karun.ipm.ac.ir \\
Fax: 009821-2290648
\end{center}

\vspace*{1cm}

\begin{abstract}
In this paper, we give a complete analysis of an SIS
epidemiological model in a population of varying size with two
dissimilar groups of infective individuals. It is mainly based on
the discussion of the existence and stability of equilibria of
the proportions system and the result is in terms of a threshold
parameter which governs the stability of the disease free
equilibrium.
\end{abstract}


\section{Introduction}
The social mixing structure of a population or a group of
interacting populations play a crucial role in the dynamics of a
disease transmission. (See \cite{3} and references therein.) In
almost all attempts to combine epidemiological data with
mathematical modeling, there has been a recognition of the need
to consider the structure of social interactions among the
individuals in the populations. (See \cite{11} and references
therein.) Many authors have considered the multigroup models in
which heterogenus subpopulations may participate to the epidemic
process with different parameters \cite{2}. For SIS type models,
a rather complete analysis of existence and global stability of a
nontrivial epidemic state has been carried out by Lajmanovich and
York \cite{7}. In their work, the size of each subpopulation is
assumed to be constant.

A famous example for these subpopulations is the core group, i.e.
the highly sexually active subgroups \cite{5}. It has become
increasingly clear that the transmission within and among core
subgroups is an important factor in the transmission of HIV/AIDS
\cite{6}. In order to consider the core group in an SIS
epidemiological model, we divide the population into two subgroups
each of them consists of susceptible and infective individuals.
One of these subpopulations can be viewed as the core group. In
this paper we consider another type of core group that is {\it
post-infection core group }, i.e. individuals that become part of
the core group after being infected. This hypothesis is plausible
for a contagious and fatal disease like AIDS. From the
psychological perspective, this group may be classified as
violent.

In  this paper we examine an SIS model of disease transmission in
a population of varying size with two dissimilar groups of
infective individuals. One of these groups can be viewed as the
post-infection core group . We also assume that the birth rate of
susceptibles may be more than that of infectives. This is similar
to the demographic assumption in \cite{8}. This paper is mainly
based on the discussion of the existence and stability of
equilibria of the proportions system. First of all, in the next
section, we introduce the model and some concepts of ODE's
related to the system. In Section 3, we present some basic
results concerning the nonexistence of certain types of
solutions. In Section 4, we give a complete global analysis of
the proportions system which is reduced to a planar system. The
result is in terms of a threshold parameter which governs the
stability of the disease free equilibrium.

\section{The Model}
In order to derive our model, we divide the population into three
groups: Susceptibles, $S$, and two groups of infectives, $I_1$
and $I_2$. We set $N=S+I_1+I_2$ which is the total size of the
population and we use the following parameters which are assumed
to be positive unless otherwise specified:

$b$: per capita birth rate of susceptibles,

$b_1$: per capita birth rate of infectives which is assumed to be $\leq b$,

$d$: per capita disease free death rate,

$\varepsilon$: excess per capita death  rate of infectives,

$\lambda_1$: effective per capita contact rate of $I_1$,

$\lambda_2$: effective per capita contact rate of $I_2$,

$\gamma_1$: per capita recovery rate of $I_1$,

$\gamma_2$: per capita recovery rate of $I_2$. \\
We also assume that the susceptible individuals which have been infected,
enter to the group $I_1$ and $I_2$ of proportions $p$ and $q$ respectively,
hence $p+q=1$.

The above hypotheses leads to the following system of
differential equations in $\R^3_+$, where ``$'$''
 denotes the derivatives with respect to $t$, the time,
$$\left \{ \begin{array}{lll}
{S'}=b_1N+(b_2-d)S+\gamma_1 I_1+\gamma_2 I_2-\lambda_1
\frac{I_1S}{N}-\lambda_2 \frac{I_2S}{N}, &\hfill{} & (2-1)\\
{I'}_1=p\left (\lambda_1 \frac{I_1S}{N} +\lambda_2 \frac{I_2S}{N}
\right)-(d+\varepsilon +\gamma_1) I_1, &\hfill{} & (2-2)\\
{I'}_2=q\left (\lambda_1 \frac{I_1S}{N} +\lambda_2 \frac{I_2S}{N}
\right)-(d+\varepsilon +\gamma_2) I_2, &\hfill{} & (2-3)
\end{array}\right .$$
where $b_2=b-b_1$ and $\frac{\lambda_j I_j}{N}$ is of the
proportionate (or random) mixing type \cite{9}, \cite{4}. By
adding the above three equations, the total population equation is
$${N'}=(b_1-d)N+b_2 S-\varepsilon (I_1+I_2)$$
Setting $s=\frac{S}{N}, i_1=\frac{I_1}{N}$ and $i_2=\frac{I_2}{N}$, we arrive
at the following system of equations:
$$\left \{\begin{array}{lll}
 s'=b_1(1-s)+b_2 s(1-s)+\gamma_1 i_1+\gamma_2
i_2+(\var-\lambda_1)i_1 s+ (\var -\lambda_2) i_2 s, &
\hfill{} & (2-1)'\\
i'_1=ps(\lambda_1 i_1+\lambda_2i_2)+\var i_1
(i_1+i_2)-(b_1+\var +\gamma_1)i_1-b_2 si_1, & \hfill{}
& (2-2)'\\
i'_2=qs (\lambda_1 i_1+\lambda_2i_2)+\var
i_2(i_1+i_2)-(b_1+\var+\gamma_2)i_2-b_2 si_2. &
\hfill{} & (2-3)'
\end{array}\right .$$

In order to determine
 the asymptotic behaviour of the solutions of this system of
equations, we need the following concepts of $ODE's$ related
to our system.

Given an autonomous system of ordinary differential equations in $\R^n$,
\begin{equation}
\frac{dx}{dt}=f(x),
\tag{$2-5$}
\end{equation}
we will denote by x.t the value of the solution of this system at time $t$,
that is $x$ initially.
For $V\subseteq \R^n$, $J\subseteq \R$, we let $V.J=\{x.t: x\in V, t\in J\}$.
The set $V$ is called positively invariant if $V.\R^+=V$.
 For $Y\subseteq \R^n$
the $\omega$-limit (resp. the $\alpha$-limit) set of $Y$ is
defined to be the maximal invariant set in the closure of $Y. [0,
\infty)$ (resp. $Y.(-\infty,0]$). We say that $\gamma(t)$ is an
orbit running from $x_0$ to $x_1$ if $\D{\lim_{t\rg-\infty}}
\gamma(t)=x_0$ and $\D{\lim_{t\rg+\infty}} \gamma(t)=x_1$. These
two points must be equilibria and such an orbit is called
heteroclinic orbit. When $x_1$ coincides with $x_0$, it is called
a homoclinic orbit. A closed curve connecting several equilibria
whose segments between successive equilibria are heteroclinic
orbits is called a phase polygon. By a sink we mean an
equilibrium at which all the eigenvalues of the linearized system
have negative real parts. Such a point is called a source if all
of these eigenvalues have positive real parts. If some of these
eigenvalues have positive real parts and the others negative real
parts, then the equilibrium
 is called a saddle point and it is called nondegenerate if all
of these eigenvalues are nonzero.

\section{Some Basic Results}
We start our analysis with some basic results about the system
$(2-1)'-(2-3)'$. If we set $\Sigma = s+{i}_1+{i}_2 $, then
$\Sigma'=(1-\Sigma)(b_1+b_2 s-\var i_1 -\var i_2)$. Therefore the
plane $\sum=1$ is invariant. We consider the feasibility region
$$D=\{(s, i_1, i_2): s+ i_1+i_2=1, s\geq 0 ,
i_1\geq 0, i_2\geq 0 \}$$
which is a triangle and on its sides we have:
$$\left \{ \begin{array}{l}
s=0 \Longrightarrow s'=b_1+\gamma_1 i_1+\gamma_2 i_2,\\
i_1=0 \Longrightarrow i'_1=p\lambda_2 si_2,\\
i_2=0\Longrightarrow i'_2=q\lambda_1 si_1.
 \end{array}\right.$$
It follows that $D$ is positively invariant and the disease free
equilibrium $(1,0,0)$ is the only rest point  on $\partial D$,
the boundary of $D$. Indeed our vector field points inward on
$\partial D -\{(1,0,0)\}$. So every solution of the system
$(2-1)'-(2-3)'$ which starts in $\partial D-\{(1,0,0)\}$,
immediately gets into $\st{\circ}{D}$, the interior of $D$.

From now on, we examine the dynamics of this system in the
feasibility region $D$. The following theorem is a modification
of Theorem 4.1 in \cite{1}, concerning the nonexistence of certain
types of solutions.
\paragraph{Theorem 3.1.} Let $f$ be a smooth vector field in $\R^3$ and
$\gamma(t)$ be a closed piecewise smooth curve which is the
boundary of an orientable smooth surface $S\subset \R^3$. Suppose
$g:U\rg \R^3$ is defined and is smooth in a neighborhood $U$ of
$S$. Moreover it satisfies $g(\gamma(t)).f(\gamma(t))\geq 0$ and
(curl g) $\cdot n <0$, where $n$ is  the unit normal to $S$. Then
$\gamma$ is not a finite union of the orbits of the system (2-5).

In order to apply the above theorem, we define $g=g_1+g_2+g_3$ by
$$g_1 (i_1, i_2)= \left [ 0, - \frac{ f_3 (i_1,
i_2)}{i_1 i_2} , \frac{f_2(i_1, i_2)}{i_1
i_2} \right], $$
$$g_2 (s, i_2)= \left [\frac{ f_3 (s,
i_2)}{si_2}
 , 0,  -\frac{f_1(s, i_2)}{s
i_2} \right], $$
$$g_3 (s, i_1)= \left [-\frac{ f_2 (s,
i_1)}{si_1}
 ,  -\frac{f_1(s, i_1)}{s
i_1},0  \right], $$
where $f_1,f_2$ and $f_3$ deduced by $\Sigma=1$ on
 the right hand side of $(2-1)', (2-2)'$ and
$(2-3)'$ respectively. Now after  some computations we get
$$(\text{curl} \ g).(1,1,1)=- \left ( \frac{p \lambda_2}{i_1^2}+ \frac{q
\lambda_1}{i_2^2}+\frac{b_1+\gamma_1}{i_2
s^2}+\frac{b_1+\gamma_2}{i_1 s^2} \right).$$
\paragraph{Corollary 3.2.} The system $(2-1)'-(2-3)'$ has no periodic orbits,
homoclinic orbits or phase polygons in $\st{\circ}{D}$.
\paragraph{Proof.} We use Theorem 3.1. for $f=(f_1,f_2, f_3)$. Here we have
$g.f=0$ and $(\text{curl}\ g) .(1,1,1)<0$ in $\st{\circ}{D}$.
$\square$
\paragraph{Lemma 3.3.} The $\omega$-limit set of each orbit of the system
$(2-1)'-(2-3)'$ with initial point in $D$ is a rest point.
\paragraph{Proof.} Suppose the contrary, then the $\omega$-limit set has
a regular point in
$\st{\circ}{D}$. Let $x$ be such a point and $h$
 be its first return map. For
a point $y$ near $x$ on the transversal, let $V$ be the region surrounded by
the orbit $\gamma$ from $y$ to $h(y)$ and the segment between them.
This region is known as Bendixon sack. (See Figure 3.1.)

Now by Stokes' theorem
$$\int\int_V
 (\text{curl}\ g).(1,1,1)d\sigma =\int_\gamma g.f dt+\int_0^1 g(ty+
(1-t)h(y)).(y-h(y))dt.$$
Since $g.f=0$ and $h(x)=x$, the right hand side of the above equality tends to
zero when $y$ tends to $x$. But the left hand side tends to the integral over
the region bounded by the $\omega$-limit set. This is a contradiction since
$(\text{curl}\ g).(1,1,1)<0$ in $\st{\circ}{D}$.~$\square$\\

{\large {
\unitlength 1.00mm
\linethickness{1.1pt}
\begin{picture}(77.33,80.67)
\put(35.33,80.67){\line(1,-1){17.33}}
\bezier{196}(68.67,63.67)(77.33,44.33)(50.00,38.00)
\bezier{100}(28.33,57.33)(27.67,67.67)(37.67,78.00)
\bezier{156}(50.00,38.00)(31.00,37.33)(28.33,56.67)
\put(35.3,75.00){\vector(3,4){1.33}}
\bezier{136}(68.67,63.67)(59.67,79.33)(46.67,69.67)
\put(49.33,71.12){\vector(4,3){1.67}}
\end{picture}}}

\vspace*{-3cm}

\begin{center}
Fig. 3.1. The Bendixon sack.
\end{center}

\paragraph{Remark 3.4.} When the $\omega$-limit set lies in $\st{\circ}{D}$
the above result is easily concluded by the generalized
Poincar$\acute{\text{e}}$-Bendixon theorem \cite{10} and Corollary
4.2. Similarly if the $\alpha$-limit set of an orbit of the system
$(2-1)'-(2-3)'$ lies in
$\st{\circ}{D}$, it must be a single point. \\

\section{The Planar System}
Using the equality $s+i_1+i_2=1$, we see that our system is essentially two
dimensional. Thus we can eliminate one of the variables, say $s$,
 to arrive at the
following quadratic planar system
$$\left \{ \begin{array}{lll}
i'_1=(p\lambda_1-b-\var -\gamma_1)i_1+p\lambda_2
i_2+(i_1+i_2)((b_2+\var-p\lambda_1)i_1-p\lambda_2 i_2), & \hfill{} & (3-1)\\
i'_2=q\lambda_1
i_1+(q\lambda_2-b-\var-\gamma_2)i_2+(i_1+i_2)((b_2+\var-q\lambda_2)i_2-q
\lambda_1 i_1).  &\hfill{} &(3-2)
\end{array}\right . $$
The dynamics of the system $(2-1)'-(2-3)'$ on $D$ is equivalent to the
dynamics of
 this planar system in the positively invariant region
$$D_1=\{(i_1,i_2)|i_1\geq 0, i_2\geq 0, i_1+i_2\leq 1 \}.$$
The matrix of the linearization of the system $(3-1),(3-2)$ at the origin is:
$$C=\left [ \begin{array}{lc}
p \lambda_1-b-\var-\gamma_1 & p \lambda_2\\
q\lambda_1 & q\lambda_2-b-\var -\gamma_2
\end{array} \right ],$$
with $det\ C=
(b+\var+\gamma_1)(b+\var+\gamma_2)-p\lambda_1(b+\var+\gamma_2)
-q\lambda_1(b+\var+\gamma_1)$. We set
$R_0=\frac{p\lambda_1}{b+\var+\gamma_1}+\frac{q\lambda_2}{b+\var+\gamma_2}$.
Hence if $R_0<1$,
 then $det\ C>0$ and $trace\ C<0$ and if $R_0>1$ then $det\ C<0$.
Thus we have proved the following lemma.
\paragraph{Lemma 4.1.} Let $R_0$ be the above threshold.
Then the origin is a sink (resp. a saddle) for the system
$(3-1),(3-2)$ whenever $R_0<1$ (resp. $R_0>1$).
\paragraph{Lemma 4.2.} The trace of the linearization of the system
$(3-1),(3-2)$ at a rest point in $\st{\circ}{D}_1$ is negative.
\paragraph{Proof.} We compute the trace at a rest point in $\st{\circ}{D}_1$.
$$\begin{array}{l}
\frac{\partial i'_1}{\partial
i_1}=p\lambda_1-b-\var-\gamma_1+(b_2+\var-p\lambda_1 -p\lambda_2)
i_2+2(b_2+\var-p \lambda_1)i_1,\\
\frac{\partial i'_2}{\partial
i_2}=p\lambda_2-b-\var-\gamma_2+(b_2+\var-q\lambda_2 -q\lambda_1)
i_1+2(b_2+\var-q \lambda_2)i_2.
\end{array}$$
From $i'_1=0$ and $i'_2=0$, we get
$$\begin{array}{l}
\frac{\partial i'_1}{\partial
i_1}=-p\lambda_2\frac{i_2}{i_1}+p\lambda_2\frac{i_2^2}{i_1}+
(b_2+\var-p\lambda_1)i_1=-p\lambda_2 \frac{i_2}{i_1}
(1-i_2)+(b_2+\var-p\lambda_1)i_1,\\
\frac{\partial i'_2}{\partial
i_2}=-q\lambda_1\frac{i_1}{i_2}+q\lambda_1\frac{i_1^2}{i_2}+
(b_2+\var-q\lambda_2)i_2=-q\lambda_1 \frac{i_1}{i_2}
(1-i_1)+(b_2+\var-q\lambda_2)i_2.
\end{array}$$
Using the equality $s+i_1+i_2=1$, we obtain
$$\frac{\partial i'_1}{\partial i_1}+\frac{\partial i'_2}{\partial
i_2}=(b_2+\var-\lambda_1)i_1+(b_2+\var-\lambda_2)i_2-p\lambda_2
\frac{i_2s}{i_1} - q\lambda_1 \frac{i_1s}{i_2}.$$ Now from
$(2-1)'$ we have
$$s'=b_1(i_1+i_2)+b_2s(i_1+i_2)+\gamma_1i_1+\gamma_2i_2+(\var-
\lambda_1)i_1s+(\var-\lambda_2)i_2s=0.$$ Thus
$(b_2+\var-\lambda_1)i_1s+(b_2+\var-\lambda_2)i_2 s<0$ and it
follows that $\frac{\partial i'_1}{\partial i_1}+\frac{\partial
i' _2}{\partial i_2}<0$. $\square$

The following two corollaries are immediate results of the above lemma.
\paragraph{Corollary 4.3.} The system $(3-1),(3-2)$ has no
source in $\st{\circ}{D}_1$.
\paragraph{Corollary 4.4.} Every nondegenerate
rest point of  the system $(3-1),(3-2)$ in $\st{\circ}{D}_1$ is
hyperbolic.
\paragraph{Remark 4.5.} A nondegenerate rest point of the system
$(3-1),(3-2)$ is obtained by a transversal intersection of the
two conic sections $i'_1=0$ and $i'_2=0$.
\paragraph{Proposition 4.6.} There is at most one rest point in
$\st{\circ}{D}_1$ for the system $(3-1),(3-2)$. Moreover such a
rest point is always hyperbolic.
\paragraph{Proof:} From the equilibrium conditions $i'_1=i'_2=0$, we get the
following equation which is homogeneous with respect to $i_1$ and
$i_2$ of second order.
$$((p\lambda_1-b-\var-\gamma_1)i_1+p\lambda_2 i_2)((b_2+\var
-q\lambda_2)i_2-q\lambda_1 i_1)=$$
$$(q\lambda_1i_1+(q\lambda_2-b-\var-\gamma_2)i_2)((b_2+\var-p\lambda_1)i_1-p
\lambda_2 i_2).$$ This equality can be written as
$$q\lambda_1(b_1+\gamma_1)i_1^2+(*) i_1i_2-p \lambda_2 (b_1+\gamma_2)i_2^2=0$$
where $(*)$ is a statement in terms of the involved parameters.
The set of all roots of this quadratic equation consists of two
lines through the origin in the $(i_1,i_2)$ plane. One of these
lines has negative slope and meets $D_1$ only at the origin. Thus
the other line contains all rest points of the system
$(3-1),(3-2)$ in $D_1$. Since each line contains at most two rest
points of a quadratic planar system and this line contains the
origin, it follows that $\st{\circ}{D}_1$ contains at most one
rest point. This rest point is obtained by a transversal
intersection of this line and each of
 the conic sections $i'_1=0$ or $i'_2=0$. It is easy to see that
at this rest point, these two conic sections must intersect
 transversally. Now
by Remark 4.5, this rest point is nondegenerate and by Corollary 4.4, it
must be hyperbolic. $\square$
\paragraph{Remark 4.7.} We have indeed shown that all rest points of the
system $(3-1),(3-2)$ which are not more than three points, are nondegenerate,
except the origin in the case $R_0=1$.

Now we are ready to prove our main result about the dynamics of the system
$(2-1)'-(2-3)'$ in $D$.
\paragraph{Theorem 4.8.} (i) If  $R_0\leq 1$, then $(1,0,0)$ is a
global attractor in $D$ \\ (ii) If $R_0>1$, then there exists a
unique rest point (an endemic equilibrium) in $\st{\circ}{D}$
which attracts $D-\{(1,0,0)\}$.
\paragraph{Proof.} When $R_0<1$, the origin is a sink for the planar system
$(3-1),(3-2)$. If there exists another rest point in
$\st{\circ}{D}_1$ for this system, it must be unique and
hyperbolic. By Corollary 4.3 it cannot be a source. If it is a
sink, then we will have two sinks in $D_1$. The basins of
attraction of these two points are open and by Lemma 3.3, $D_1$
is the union of these two open subsets. This contradicts the
connectedness of $D_1$ and shows that it cannot be a sink. Now
suppose that there is a saddle point in $\st{\circ}{D}_1$. By
Corollary 3.2, there is no homoclinic orbit in $\st{\circ}{D}_1$.
Hence the origin attracts the unstable manifold of the saddle
point.  Now the region bounded by the unstable manifold  contains
some part of the stable manifold of of the saddle point. Thus the
$\alpha$-limit set of this part of the stable manifold is a rest
point in $D_1$ by Remark 3.4. This is a contradiction since this
rest point can be neither the origin nor the saddle point.
 Therefore the origin is the only
rest point in $D_1$ for the system $(3-1),(3-2)$ and by Lemma
2.3, it is the $\omega$-limit set of all points of $D_1$.

The above fact is still valid for the limiting  case, $R_0=1$. To
see this, suppose that there exists another rest point in $D_1$,
then it must be hyperbolic and belong to $\st{\circ}{D}_1$. Thus
it remains in $\st{\circ}{D}_1$ when the involved parameters are
slightly changed to get $R_0<1$ which contradicts the above
result. This finishes the proof of (i).

Now suppose $R_0>1$. Then the origin is a saddle point for the
planar system $(3-1),(3-2)$. Thus by Lemma 2.3, there must be some
rest point in $\st{\circ}{D}_1$. Since such a rest point is
unique and hyperbolic, it must be a sink and attract all points
of $D_1$ except the stable manifold of the origin. We claim that
the stable manifold meets $D_1$ only at the origin. To see this
notice that some part of the unstable manifold of the origin must
be outside of $D_1$ for its right angle. Since $D_1$ is positively
invariant, the stable manifold does not intersect
$\st{\circ}{D}_1$. Moreover, the vector field points inward on
$\partial D_1-\{(0,0)\}$. Thus the stable manifold of the origin
does not intersect $\partial D_1-\{(0,0)\}$ either. This shows
that the origin cannot attract any point of $D_1-\{(0,0)\}$. Thus
the unique sink attracts $\st{\circ}{D}_1-\{(0,0)\}$. It means
that there is a unique rest point in $\st{\circ}{D}$ for the
system $(2-1)'-(2-3)'$ which attracts $D-\{(1,0,0)\}$. $\square$

\paragraph{Remark 4.9.} In the above argument, in
order to prove the global asymptotic stability of the endemic
equilibrium (i.e. the unique rest point in $\st{\circ}{D}_1$), we
showed that the stable manifold of the origin cannot intersect
$\st{\circ}{D}_1$. It is a special case of the following fact.
Let $X$ be a smooth vector field on a smooth manifold $M$ and
$D\subset M$ is a positively invariant region with a saddle point
on $\partial D$. If the unstable manifold of this saddle point
contains a point of $(M-D)^{\circ}$, then its stable manifold
cannot intersect $\st{\circ}{D}$. In order to prove it, one can
follow our proof in the special case and observe that this is a
direct consequence of the Hartman-Grobman theorem. However, it is
obvious by the Inclination Lemma \cite{10}.

\paragraph{Remark 4.10.} If we consider $I_2$ as the post-infection
core group, then the threshold $R_0=\frac{p\lambda_1}
{b+\var+\gamma_1}+\frac{q\lambda_2}{b+\var+\gamma_2}$ clearly
shows the effect of this group on the epidemics process. Although
the probability $q$ is a small number, the fact $\lambda
_2>>\lambda _1$ causes the term $q\lambda _2$ to be significant.

\paragraph{Acknowledgments.} The author  would like
to thank Institute for Studies in Theoretical Physics and
Mathematics for supporting this research.

\end{document}